\documentclass{amsart}
\usepackage{amssymb}
\usepackage{a4}
\usepackage{color}

\newtheorem{theorem}{Theorem}[section]
\newtheorem{lemma}[theorem]{Lemma}
\newtheorem{remark}[theorem]{Remark}

\newtheorem{example}[theorem]{Example}
\newtheorem{de}[theorem]{Definition}

\newcommand{\g}[2]{\ensuremath{\langle #1,#2 \rangle}}
\newcommand{\jxa}{\ensuremath{\mathcal{J}(\pi_x)}}
\newcommand{\jya}{\ensuremath{\mathcal{J}(\pi_y)}}
\begin{document}
\title[Complex Osserman K{\"a}hler manifolds]
 {Complex Osserman K{\"a}hler manifolds in dimension four}
\author{M. Brozos-V\'azquez and P. Gilkey}
\address{MBV: E. U. P. Ferrol, Mathematics Department, University of A Coru\~na, Spain\\
E-mail: mbrozos@udc.es}
\address{PG: Mathematics Department, University of Oregon\\
  Eugene OR 97403 USA\\
  E-mail: gilkey@uoregon.edu}

\begin{abstract}
Let $\mathcal{H}$ be a $4$-dimensional almost-Hermitian manifold which satisfies the K\"ahler identity. We show that
$\mathcal{H}$ is complex Osserman if and only if $\mathcal{H}$ has constant holomorphic sectional curvature. We also
classify in arbitrary dimensions all the complex Osserman K\"ahler  models which do not have 3
eigenvalues.
\end{abstract}
\maketitle

\section{Introduction}

\subsection{The geometric context} Let $R(x,y):=\nabla_x\nabla_y-\nabla_y\nabla_x-\nabla_{[x,y]}$ be the curvature operator of
a Riemannian manifold
$\mathcal{M}:=(M,g)$ of dimension
$n$ and let $R(x,y,z,w):=g(R(x,y)z,w)$ be the associated curvature tensor. The eigenvalue structure of
various operators naturally associated to
$R$ has been studied in the last decade to obtain geometric information about
$\mathcal{M}$. We recall some notation:
\begin{de}\label{defn-1.1}
Let $\mathcal{M}=(M,g)$ be a Riemannian manifold.
\begin{enumerate}
\item If $R(x,y,y,x)=c$ for any orthonormal set $\{x,y\}$, then $R$ is said to have {\rm constant sectional curvature}.
\item The {\rm Jacobi operator} is defined by $\mathcal{J}(x):y\rightarrow R(y,x)x$.
\item Let $\{e_1,...,e_k\}$ be an orthonormal basis for a subspace $\sigma$. Following
\cite{StanilovVidev}, the {\rm higher order Jacobi operator} $\mathcal{J}(\sigma)$ is defined by setting:
$$\mathcal{J}(\sigma):=\mathcal{J}(e_1)+...+\mathcal{J}(e_k)\,.$$
This is independent of the particular orthonormal basis chosen for $\sigma$.
\item
$\mathcal{M}$ is said to be {\rm Osserman} if the eigenvalues of $\mathcal{J}(x)$ are constant
on the bundle of unit tangent vectors $S(\mathcal{M})$.
\end{enumerate}\end{de}
Any  locally $2$-point homogeneous space is clearly Osserman. Osserman \cite{Osserman} wondered if the
converse held. This converse implication has been established except for the (possibly) exceptional case
$n=16$; the Jacobi operator of $\mathcal{M}$ has constant eigenvalues on $S(\mathcal{M})$ if and only if
$\mathcal{M}$ is locally a 2-point homogeneous space
\cite{Chi,Gilkey-Swann-Vanhecke,nikolayevsky-1,nikolayevsky-2}. There are related problems defined by other
natural operators. The conformal Jacobi operator has been investigated \cite{b-g-04,b-g-05,nikolayevsky-3}, the
skew-symmetric curvature operator has been investigated
\cite{GLS,IP-98,N-04}, and the higher order Jacobi operator has been investigated \cite{gilkey-p-Osserman}; we
refer to
\cite{GKV,Gilkey-book} for further details.

We now pass to the complex setting. Recall that $\mathcal{H}:=(M,g,J)$ is said to be an {\it almost-Hermitian
manifold} if $J$ is an endomorphism of the tangent bundle $TM$ which satisfies $J^2=-\operatorname{Id}$ and
$J^*g=g$.  Again, we recall some notation:

\begin{de}\label{defn-1.2}
Let $\mathcal{H}$ be an almost-Hermitian manifold.\begin{enumerate}
\item A $2$-dimensional subspace $\pi$ of $TM$ is said to be a {\rm complex line} if $J\pi=\pi$. Let
$\mathbb{CP}(\mathcal{H})$ be the bundle of complex lines. The map $x\rightarrow\pi_x:=\operatorname{Span}\{x,Jx\}$
defines the {\rm Hopf fibration} $S(\mathcal{H})\rightarrow\mathbb{CP}(\mathcal{H})$. If $\pi\in\mathbb{CP}(\mathcal{H})$, the
{\rm complex curvature operator} and the {\rm complex Jacobi operator} are defined, respectively, by setting $R(\pi):=R(x,Jx)$
and
$\mathcal{J}(\pi):=\mathcal{J}(x)+\mathcal{J}(Jx)$; these operators are independent of the particular unit vector
$x\in S(\pi)$ which is chosen.
\item Let
$s(x):=R(x,Jx,Jx,x)$ be the holomorphic sectional curvature; $\mathcal{H}$ is
said to have {\rm constant holomorphic sectional curvature} if $s(\cdot)=c$ on $S(\mathcal{H})$.
\item If
$R(Jx,Jy)=R(x,y)$ for all $x,y$, then $\mathcal{H}$ is said to {\rm satisfy the K{\"a}hler identity}.
\item An
endomorphism $\Xi$ of $TM$ is said to be {\rm complex} if $J\Xi=\Xi J$.
\item $\mathcal{H}$ is said to be {\rm complex Osserman} if
 $\mathcal{J}(\pi)$ is complex for every $\pi\in\mathbb{CP}(\mathcal{H})$ and if the eigenvalues of $\mathcal{J}(\pi)$ are
constant on $\mathbb{CP}(\mathcal{H})$.
In this setting, the eigenvalues and eigenvalue multiplicities of $\mathcal{J}(\pi)$ for any
$\pi\in\mathbb{CP}(\mathcal{H})$ are said to be {\rm the eigenvalues and eigenvalue multiplicities} of $\mathcal{H}$.
\end{enumerate}\end{de}

We shall be assuming for the most part that $\mathcal{H}$ satisfies the K\"ahler identity; Lemma \ref{lem-3.1} below shows
that necessarily
$\mathcal{J}(\pi)$ is complex in this setting.

\subsection{The algebraic context} Let $V$ be an $n$-dimensional real vector space. An element
$A\in\otimes^4V^*$ is said to be an {\it algebraic curvature tensor} if $A$ has the symmetries of
the Riemann curvature tensor, i.e. if for all $x,y,z,w\in V$, we have:
\begin{eqnarray}
&&A(x,y,z,w)=-A(y,x,z,w)=A(z,w,x,y)\,,\label{eqn-1}\\
&&A(x,y,z,w)+A(y,z,x,w)+A(z,x,y,w)=0\,.\label{eqn-2}
\end{eqnarray}
One says that $\mathfrak{H}:=(V,\langle\cdot,\cdot\rangle,J,A)$ is an {\it almost-Hermitian curvature model} if
$A$ is an algebraic curvature tensor on $V$, if $\langle\cdot,\cdot\rangle$ is a positive definite inner product on $V$, and
if $J$ is an endomorphism of $V$ satisfying
$J^2=-\operatorname{Id}$ and
$J^*\langle\cdot,\cdot\rangle=\langle\cdot,\cdot\rangle$.
Let $A(x,y)$ be the
corresponding curvature operator. The notions of Definition \ref{defn-1.1} and Definition \ref{defn-1.2} then extend
immediately to this setting.

An almost-Hermitian manifold $\mathcal{H}$ is said to be K{\"a}hler if $\nabla J=0$; such a manifold satisfies the
{\it K{\"a}hler identity} $R(Jx,Jy)=R(x,y)$ discussed above. An almost-Hermitian curvature model $\mathfrak{H}$ is said to be
a {\it K\"ahler model} if $\mathfrak{H}$ satisfies the K\"ahler identity. Every K\"ahler model can be geometrically realized
by a K{\"a}hler manifold \cite{bv-g-m}.

The eigenvalue structure of a complex Osserman K\"ahler  model is very restrictive. We shall establish the following result in
Section \ref{sect-3}:
\begin{theorem}\label{thm-1.3}
Let $\mathfrak{H}$ be a complex Osserman K\"ahler model of dimension $n\ge4$ which is not flat.
Then one of the following
holds:
\begin{enumerate}
\item There are 2 eigenvalues of multiplicities $(n-2,2)$.
\item There are 3 eigenvalues of multiplicities $(n-4,2,2)$ with $n=4k\ge8$.
\end{enumerate}\end{theorem}

\subsection{Algebraic classification} Let $V$ be a vector space of dimension $n=2m$, let $\langle\cdot,\cdot\rangle$ be a
positive definite symmetric inner product on $V$. If $\phi\in S^2(V^*)$ is a symmetric 2-tensor, we define
$$A_\phi(x,y,z,w)=\phi(x,w)\phi(y,z)-\phi(x,z)\phi(y,w)\,.$$
Similarly if $\psi\in\Lambda^2(V^*)$ is an anti-symmetric $2$-tensor, we define:
$$A_\psi(x,y,z,w)=\psi(x,w) \psi(y,z)-\psi(x,z)\psi(y,w)-2\psi(x,y)\psi(z,w)\,.$$
It is an easy calculation to show that $A_\phi$ and $A_\psi$ are algebraic curvature tensors -- see, for example, Lemma 1.8.1 of
\cite{Gilkey-book}. Let $J$ be an endomorphism of $V$ with
$J^2=-\operatorname{Id}$ and with
$J^*\langle\cdot,\cdot\rangle=\langle\cdot,\cdot\rangle$. The following examples will play a crucial role in our
development.
\begin{example}
\rm Let $\mathfrak{A}_\mu^n=(V,\langle\cdot,\cdot\rangle,\mu A_0)$ be an $n$-dimensional real model where $A_0=A_\phi$ is defined by taking
$\phi=\langle\cdot,\cdot\rangle$, i.e.:
$$
A_0(x,y,z,w):= \langle x,w\rangle \langle y,z\rangle-\langle x,z\rangle \langle y,w\rangle\,.
$$
\end{example}

\begin{example}\label{exm-1.5}
\rm Let $\mathfrak{B}_\mu^n=(V,\langle\cdot,\cdot\rangle,J,\frac14\mu(A_0+A_J))$ where $A_J$ is defined by taking $\psi:=\langle
J\cdot,\cdot\rangle$ to be negative of the K\"ahler form, i.e.:
$$
A_J(x,y,z,w):= \langle Jx,w\rangle \langle Jy,z\rangle-\langle Jx,z\rangle \langle Jy,w\rangle
    -2\langle Jx,y\rangle \langle Jz,w\rangle\,.
$$
The non-zero curvatures are determined  for $i<j$, up to the usual $\mathbb{Z}_2$ symmetries, by the relations:
$$
\begin{array}{ll}
B_\mu^n(x_i,Jx_i,Jx_i,x_i)=\mu,&B_\mu^n(x_i,x_j,x_j,x_i)=\frac14\mu,\\
B_\mu^n(Jx_i,x_j,x_j,Jx_i)=\frac14\mu,& B_\mu^n(Jx_i,Jx_j,Jx_j,Jx_i)=\frac14\mu,\vphantom{\vrule height 11pt}\\
B_\mu^n(x_i,x_j,Jx_j,Jx_i)=\frac14\mu,&B_\mu^n(x_i,Jx_j,Jx_i,x_j)=\frac14\mu,\vphantom{\vrule height 11pt}\\
B_\mu^n(x_i,Jx_i,Jx_j,x_j)=\frac12\mu\,.\vphantom{\vrule height 11pt}
\end{array}$$
\end{example}
\begin{example}\label{exm-1.6}
\rm Fix an orthonormal basis $\{x_1,Jx_1,...,x_m,Jx_m\}$ on the vector space $V$ of dimension $n$. Let
$\mathfrak{C}_\mu^n=(V,\langle\cdot,\cdot\rangle,J,C_\mu^n)$ where the non-zero components of $C_\mu^n$ are determined by the following
identities modulo the usual $\mathbb{Z}_2$ symmetries for
$i\ne j$:
$$
\begin{array}{ll}
C_\mu^n(x_i,Jx_i,Jx_i,x_i)=\mu,&C_\mu^n(x_i,x_j,x_j,x_i)=\frac{-\mu}{2},\\
C_\mu^n(Jx_i,x_j,x_j,Jx_i)=\frac{\mu}{2},&C_\mu^n(Jx_i,Jx_j,Jx_j,Jx_i)=\frac{-\mu}{2},\vphantom{\vrule height 11pt}\\
C_\mu^n(x_i,x_j,Jx_j,Jx_i)=\frac{-\mu}{2},&C_\mu^n(x_i,Jx_j,Jx_i,x_j)=\frac{\mu}{2},\vphantom{\vrule height 11pt}\\
C_\mu^n(x_i,Jx_i,Jx_j,x_j)=0\,.\vphantom{\vrule height 11pt}
\end{array}$$\end{example}

We will verify in Section \ref{sect-2} that $\mathfrak{B}_\mu^n$ and $\mathfrak{C}_\mu^n$ are complex Osserman K\"ahler
models; as they are flat if $\mu=0$, we shall usually restrict to the case $\mu\ne0$. The real model
$\mathfrak{A}_\mu^n$ has constant sectional curvature $\mu$. The complex model
$\mathfrak{B}_\mu^n$ has constant holomorphic sectional curvature $\mu$ and satisfies the K\"ahler identity. One has the
following converse -- see, for example,
\cite{KN,n73}:
\begin{lemma}\label{lem-1.7}
\ \begin{enumerate}
\item If $\mathfrak{R}=(V,\langle\cdot,\cdot\rangle,A)$ is a real model of dimension $n$ with constant sectional curvature $\mu$, then
$\mathfrak{R}$ is isomorphic to $\mathfrak{A}_\mu^n$.
\item If $\mathfrak{H}$ is a K{\"a}hler model of dimension $n$ with constant holomorphic sectional curvature $\mu$, then
$\mathfrak{H}$ is isomorphic to $\mathfrak{B}_\mu^n$.
\end{enumerate}
\end{lemma}

We will establish the following classification result in Section \ref{sect-4}:
\begin{theorem}\label{thm-1.8}
Let $n\ge4$. If $\mathfrak{H}$ is a complex Osserman K\"ahler  model of dimension $n$ with 2 eigenvalues which does not have constant
holomorphic sectional curvature, then
$\mathfrak{H}$ is isomorphic to
$\mathfrak{C}_\mu^n$ for some
$\mu$.
\end{theorem}

\subsection{Geometric classification in dimension $4$}
We shall show in Section \ref{sect-5} that $\mathfrak{C}_\mu^4$ is not geometrically realizable for $\mu\ne0$; the proof uses
the self-dual and the anti-self-dual Weyl tensors and does not generalize to the higher dimensional context. Combining this
result with Theorem
\ref{thm-1.8} then yields the main result of this paper which motivated our investigations in the first instance:

\begin{theorem}\label{thm-1.9}
Let $\mathcal{H}$ be a $4$-dimensional almost-Hermitian manifold which satisfies the K{\"a}hler identity.
Then $\mathcal{H}$ is complex Osserman if and only if
$\mathcal{H}$ has constant holomorphic sectional curvature.
\end{theorem}

\begin{remark}\rm
In fact, we shall prove a bit more. In Theorem \ref{thm-1.9}, it is only necessary to assume that $\mathcal{H}$ is pointwise
complex Osserman, i.e. the eigenvalue structure is a priori permitted to vary with the point in question. The scalar curvatures are given by
$\tau_{\mathfrak{B}_\mu^4}=6\mu$ and
$\tau_{\mathfrak{C}_\mu^4}=4\mu$. Since by Lemma \ref{lem-3.2}
$\mathcal{H}$ is Einstein, the scalar curvature (and hence $\mu$) is constant so $\mathcal{H}$ is in fact globally
complex Osserman in this setting. We also note that
$\mathcal{H}$ need not be K\"ahler to satisfy the K\"ahler identity; there are flat $4$-dimensional
almost-Hermitian manifolds which are not integrable.
\end{remark}

\section{The models $\mathfrak{A}_\mu^n$, $\mathfrak{B}_\mu^n$, and $\mathfrak{C}_\mu^n$}\label{sect-2}
In this section, we establish the basic properties of these models.

\begin{de}\rm
Let $\{y_1,y_2,y_3,y_4\}$ belong to the set $\{x_1,...,x_m,Jx_1,...,Jx_m\}$. We say that an index $i$ with
$1\le i\le m$ is an {\it impacted index} if $x_i=y_a$ or $Jx_i=y_a$ for some $a$. It is then immediate
that
$B_\mu^n(y_1,y_2,y_3,y_4)=0$ and $C_\mu^n(y_1,y_2,y_3,y_4)=0$ if there are more than 2 impacted indices or if any impacted
index appears with odd multiplicity. Finally, $B_\mu^n$ and $C_\mu^n$ vanish if $J$ appears an odd number of times.
\end{de}

\subsection{The model $\mathfrak{B}_\mu^n$}

\begin{lemma}\label{lem-2.2}
The model $\mathfrak{B}_\mu^n$ has constant holomorphic sectional curvature $\mu$ and is complex Osserman K\"ahler with 2
eigenvalues
$(\frac12\mu,\mu)$ of multiplicities
$(n-2,2)$, respectively.
\end{lemma}

\begin{proof}
As noted above $A_0$ and $A_J$ are algebraic curvature tensors. It is an easy computation that $\mathfrak{B}_\mu^n$ has constant holomorphic
sectional curvature
$\mu$. We show that
$\mathfrak{B}_\mu^n$ is K\"ahler by verifying that
\medbreak\qquad
$\frac14\mu\{A_0(x,y,Jz,Jw)+A_J(x,y,Jz,Jw)\}$
\smallbreak\qquad\quad
$=\frac14\mu\{\langle x,Jw\rangle\langle y,Jz\rangle-\langle x,Jz\rangle\langle y,Jw\rangle$
\smallbreak\qquad\qquad
$+\langle Jx,Jw\rangle\langle Jy,Jz\rangle-\langle Jx,Jz\rangle\langle Jy,Jw\rangle
    -2\langle Jx,y\rangle\langle JJz,Jw\rangle\}$
\smallbreak\qquad\quad
$=\frac14\mu\{\langle Jx,w\rangle\langle Jy,z\rangle-\langle Jx,z\rangle\langle Jy,w\rangle$
\smallbreak\qquad\qquad
$+\langle x,w\rangle\langle y,z\rangle-\langle x,z\rangle\langle y,w\rangle-2\langle Jx,y\rangle\langle Jz,w\rangle\}$
\smallbreak\qquad\quad
$=\frac14\mu\{ A_0(x,y,z,w)+A_J(x,y,z,w)\}$.
\medbreak\noindent The Jacobi operators are given by:
\begin{eqnarray}
&&\mathcal{J}_{A_0}(x)y=\left\{\begin{array}{rrr}0&\text{if}&y=x\\y&\text{if}&y\perp x\end{array}\right\},\quad
\mathcal{J}_{A_J}(x)y=\left\{\begin{array}{rrr}3y&\text{if}&y=Jx\\0&\text{if}&y\perp Jx
\end{array}\right\},\label{eqn-3}\\
&&\mathcal{J}_{\mathfrak{B}_\mu^n}(\pi_x)y=\frac\mu4\left\{\begin{array}{rrr}4y&\text{if}&y\in\pi_x
   \\2y&\text{if}&y\perp\pi_x\end{array}\right\}\,.\nonumber
\end{eqnarray}
This shows that $\mathcal{J}_{B_\mu^n}(\pi_x)$ is complex Osserman with eigenvalues $(\frac12\mu,\mu)$ of
multiplicities $(n-2,2)$, respectively.
\end{proof}

\subsection{The model $\mathfrak{C}_\mu^4$}\label{sect-2.2}
There is an auxiliary complex structure $L$, which commutes with $J$, that will be important in our investigations which
is defined by:
$$
Lx_1=x_2,\quad Lx_2=-x_1,\quad LJx_1=Jx_2,\quad LJx_2=-Jx_1\,.
$$
Let $\rho$ be the Ricci operator and let $\tau$ be the scalar curvature.

\begin{lemma}\label{lem-2.3}
Adopt the notation established above.
\begin{enumerate}
\item
$A_{\mathfrak{C}_\mu^4}=\textstyle\frac{\mu}{2}A_0+\frac{\mu}{6}A_J-\frac{\mu}{3}A_L$.
\item
$\mathfrak{C}_\mu^4$ is a complex Osserman K\"ahler model with 2 eigenvalues
$(\mu,0)$ of multiplicities $(2,2)$, respectively.
\item $
\rho_{A_0}=\rho_{A_J}=\rho_{A_L}=3 \operatorname{Id}$, $\tau_{A_0}=\tau_{A_J}=\tau_{A_L}=12$, $\tau_{\mathfrak{B}_\mu^4}=6\mu$,
$\tau_{\mathfrak{C}_\mu^4}=4\mu$.
\end{enumerate}
\end{lemma}

\begin{proof} We have
\begin{eqnarray*}
&&A_0(x,y,z,w):=\langle x,w\rangle\langle y,z\rangle-\langle x,z\rangle\langle y,w\rangle\\
&&A_J(x,y,z,w):= \langle Jx,w\rangle \langle Jy,z\rangle-\langle Jx,z\rangle \langle Jy,w\rangle
    -2\langle Jx,y\rangle \langle Jz,w\rangle,\\
&&A_L(x,y,z,w):= \langle Lx,w\rangle \langle Ly,z\rangle-\langle Lx,z\rangle \langle Ly,w\rangle
    -2\langle Lx,y\rangle \langle Lz,w\rangle\,.
\end{eqnarray*}
Let $i\ne  j$. We may verify Assertion (1) by computing:
$$\begin{array}{|r|r|r|r|r|}
\noalign{\hrule}
\text{Monomial}&A_0&A_J&A_L&C_\mu\\\noalign{\hrule}
(x_i,Jx_i,Jx_i,x_i)&1&3&0&\mu\\\noalign{\hrule}
(x_i,x_j,x_j,x_i)&1&0&3&-\frac\mu2\\\noalign{\hrule}
(Jx_i,x_j,x_j,Jx_i)&1&0&0&\frac\mu2\\\noalign{\hrule}
(Jx_i,Jx_j,Jx_j,Jx_i)&1&0&3&-\frac\mu2\\\noalign{\hrule}
(x_i,x_j,Jx_j,Jx_i)&0&1&2&-\frac\mu2\\\noalign{\hrule}
(x_i,Jx_j,Jx_i,x_j)&0&1&-1&\frac\mu2\\\noalign{\hrule}
(x_i,Jx_i,Jx_j,x_j)&0&2&1&0\\\noalign{\hrule}
\end{array}$$

Since
$A_0$,
$A_J$, and
$A_L$ are algebraic curvature tensors,
$A_{\mathfrak{C}_\mu^4}$ is an algebraic curvature tensor. The K\"ahler identity is immediate from the defining relations. We
use Equation (\ref{eqn-3}) to compute:
$$\textstyle\frac\mu2\mathcal{J}_{A_0}( \pi_x)+\frac\mu6\mathcal{J}_{A_J}( \pi_x)=\operatorname{Id}\mu\,.$$
Since $L$ is a Hermitian almost complex structure commuting with $J$, $L\pi_x$ is a complex $2$-plane as well. We complete the proof of
Assertion (2) by verifying:
\begin{eqnarray*}
&&\mathcal{J}_L(x)y=\left\{\begin{array}{rll}3y&\text{if}&y=Lx\\0&\text{if}&y\perp Lx\end{array}\right\},\quad
\mathcal{J}_L(Jx)y=\left\{\begin{array}{rll}3y&\text{if}&y=LJx\\0&\text{if}&y\perp LJx\end{array}\right\},\\
&&\mathcal{J}_L(\pi_x)y=\left\{\begin{array}{rll}3y&\text{if}&y\in L\pi_x\\0&\text{if}&y\perp L\pi_x\end{array}\right\},\quad
\mathcal{J}_{\mathfrak{C}_\mu^4}(\pi_x)y
  =\left\{\begin{array}{rll}
  0&\text{if}&y\in L\pi_x\\
 \mu y&\text{if}&y\perp L\pi_x
\end{array}\right\}\,.
\end{eqnarray*}
Assertion (3) is immediate from the definitions.
\end{proof}

\subsection{The model $\mathfrak{C}_\mu^n$}

We begin by studying the group of symmetries of the model.

\begin{lemma}\label{lem-2.4} Let $O(m)$ act diagonally on $\mathbb{R}^{2m}$. If $\Theta\in O(m)$, then
$\Theta^*\mathfrak{C}_\mu^n=\mathfrak{C}_\mu^n$.
\end{lemma}
\begin{proof} If $\Theta\in O(m)$,  then $\Theta J=J\Theta$ and $\Theta^*\langle\cdot,\cdot\rangle=\langle\cdot,\cdot\rangle$. Let $C=C_\mu^n$;
this tensor is invariant under permutations of the coordinate indices. Since
$O(m)$ is generated by coordinate permutations and by rotations in the first 2 indices, it suffices to prove the lemma for the
 elements
$$\begin{array}{rr}
\Theta_\theta(x_1)=\cos\theta x_1+\sin\theta x_2,&\Theta_\theta(Jx_1)=\cos\theta Jx_1+\sin\theta Jx_2,\\
\Theta_\theta(x_2)=-\sin\theta x_1+\cos\theta x_2,&\Theta_\theta(Jx_2)=-\sin\theta Jx_1+\cos\theta Jx_2,\\
\Theta_\theta(x_i)=x_i\text{ for }i\ge 3,&\Theta_\theta(Jx_i)=Jx_i\text{ for }i\ge3\,.
\end{array}$$We compute representative terms:
\smallbreak
$\Theta_\theta^*A(x_1,Jx_1,Jx_1,x_1)=\cos^4\theta A(x_1,Jx_1,Jx_1,x_1)+\sin^4\theta A(x_2,Jx_2,Jx_2,x_2)$
\par\quad
$+\cos^2\theta\sin^2\theta\{A(x_2,Jx_2,Jx_1,x_1)+A(x_2,Jx_1,Jx_2,x_1)+A(x_2,Jx_1,Jx_1,x_2)$
\par\quad
$+A(x_1,Jx_2,Jx_2,x_1)+A(x_1,Jx_2,Jx_1,x_2)+ A(x_1,Jx_1,Jx_2,x_2)\}$
\par\quad
$=\cos^4\theta\mu+\sin^4\theta\mu+\cos^2\theta\sin^2\theta\{0+\frac\mu2+\frac\mu2+\frac\mu2+\frac\mu2+0\}=\mu$,
\smallbreak
$\Theta_\theta^*A(Jx_1,x_2,x_2,Jx_1)=\cos^4\theta A(Jx_1,x_2,x_2,Jx_1)+\sin^4\theta A(Jx_2,x_1,x_1,Jx_2)$
\par\quad
$+\sin^2\theta\cos^2\theta\{-A(Jx_2,x_1,x_2,Jx_1)-A(Jx_2,x_2,x_1,Jx_1)$
\par\qquad$+ A(Jx_2,x_2,x_2,Jx_2)+A(Jx_1,x_1,x_1,Jx_1)$
\par\qquad
$-A(Jx_1,x_1,x_2,Jx_2)-A(Jx_1,x_2,x_1,Jx_2)\}$
\par\quad
$=(\cos^4\theta+\sin^4\theta)\frac\mu2+\cos^2\theta\sin^2\theta(-\frac\mu2-0+\mu+\mu-0-\frac\mu2\}=\frac\mu2$,
\smallbreak
$\Theta_\theta^*A(Jx_1,x_3,x_3,Jx_1)=\cos^2\theta A(Jx_1,x_3,x_3,Jx_1)+\sin^2\theta A(Jx_2,x_3,x_3,Jx_2)$
\par\quad
$=\frac\mu2(\cos^2\theta+\sin^2\theta)=\frac\mu2$,
\smallbreak
$\Theta_\theta^*A(x_1,x_2,x_2,x_1)=\cos^4\theta A(x_1,x_2,x_2,x_1)+\sin^4\theta A(x_1,x_2,x_2,x_1)$
\par\quad
$-\cos^2\theta\sin^2\theta A(x_2,x_1,x_2,x_1)-\cos^2\theta\sin^2\theta A(x_1,x_2,x_1,x_2)$
\par\quad
$=-\frac\mu2\{\cos^4\theta+2\cos^2\theta\sin^2\theta+\sin^4\theta\}=-\frac\mu2$,
\smallbreak
$\Theta_\theta^*A(x_1,x_3,x_3,x_1)=\cos^2\theta A(x_1,x_3,x_3,x_1)+\sin^2\theta A(x_2,x_3,x_3,x_2)$
\smallbreak
$=-\frac\mu2\{\cos^2\theta+\sin^2\theta\}=-\frac\mu2$,
\medbreak\noindent Since $\Theta J=J\Theta$ and since $A$ is K\"ahler, necessarily
\smallbreak
$\Theta_\theta^*A(x_i,x_j,Jx_j,Jx_i)=A(\Theta x_i,\Theta x_j,\Theta Jx_j,\Theta Jx_i)=A(\Theta x_i,\Theta x_j,J\Theta x_j,J\Theta x_i)$
\par\quad
$=A(\Theta x_i,\Theta x_j,\Theta x_j,\Theta x_i)=A(x_i,x_j,x_j,x_i)=A(x_i,x_j,Jx_j,Jx_i)$,
\smallbreak
$\Theta_\theta^*A(Jx_i,Jx_j,Jx_j,Jx_i)=A(\Theta Jx_i,\Theta Jx_j,\Theta J x_j,\Theta Jx_i)$
\par\quad$=A(J\Theta x_i,J\Theta x_j,J\Theta x_j,J\Theta x_i)=A(\Theta x_i,\Theta x_j,\Theta x_j,\Theta x_i)$
\par\quad$=A(x_i,x_j,x_j,x_i)=A(Jx_i,Jx_j,Jx_j,Jx_i)$.
\end{proof}

We remark that $\cos\theta\operatorname{Id}+\sin\theta J$ also is  an isometry of $\mathfrak{C_\mu^n}$. We now show:

\begin{lemma}\label{lem-2.5}
$\mathfrak{C}_\mu^n$ is a complex Osserman K\"ahler model with 2 eigenvalues $(\mu,0)$ of
multiplicities $(2,n-2)$, respectively. Let $\mathcal{L}:=\operatorname{Span}\{x_1,...,x_m\}$ and
let $\xi\in S(V)$.
\begin{enumerate}
\item $\mathcal{J}(\pi_\xi)\xi=\mu\xi$ if and only if $\xi\in\pi_x$ for some $x\in S(\mathcal{L})$.
\item $\mathcal{J}(\pi_\xi)\xi=0$ if and only if $\xi=(x+Jy)/\sqrt2$ for $x,y\in S(\mathcal{L})$
with $x\perp y$.
\end{enumerate}
\end{lemma}

\begin{proof} It is immediate that $C_\mu^n$ satisfies the $\mathbb{Z}_2$ symmetries of Equation (\ref{eqn-1}) and that $C_\mu^n$ satisfies
the  Bianchi identity of Equation (\ref{eqn-2}). Furthermore, one sees by inspection that
$C_\mu^n$ is K\"ahler. Thus $\mathfrak{C}_\mu^n$ is a K\"ahler model.

We wish to study the eigenvalue structure of $\mathcal{J}(\pi_\xi)$ for $\xi\in S(V)$. We expand
$\xi=\sum_ia_ix_i+\sum_jb_jJx_i$. By replacing
$\xi$ by
$J\xi$ if necessary, we may assume that $\sum a_i^2\ne0$. We use Lemma \ref{lem-2.4}. By applying an appropriate element of
$O(m)$, we may assume that $a_1\ne0$ and that
$a_i=0$ for
$i\ge2$. We then apply an appropriate element of $O(m-1)$ to assume $b_j=0$ for $j\ge3$. Thus without loss of generality, we may assume that
$$\xi=a_1x_1+b_1Jx_1+b_2Jx_2\in\mathbb{R}^4=\operatorname{Span}\{x_1,x_2,Jx_1,Jx_2\}$$
where $a_1\ne0$ and $a_1^2+b_1^2+b_2^2=1$.
Let $\eta\in\{x_3,...,x_m,Jx_3,...,Jx_m\}$.
Since $C_\mu^n$ vanishes if there are more than 2 impacted indices, we have
\begin{eqnarray*}
\mathcal{J}(\pi_\xi)\eta&=&\{\mathcal{J}(a_1x_1+b_1Jx_1)+\mathcal{J}(-b_1x_1+a_1Jx_1)\}\eta\\
&+&\{\mathcal{J}(a_2x_2+b_2Jx_2)+\mathcal{J}(-b_2x_2+a_2Jx_2)\}\eta\\
&=&(a_1^2+b_1^2)\mathcal{J}(\pi_{x_1})\eta+(a_2^2+b_2^2)\mathcal{J}(\pi_{x_2})\eta\,.
\end{eqnarray*}
One computes directly that $\mathcal{J}(\pi_{x_i})x_j=0$ for $i\ne j$. This shows $\mathcal{J}(\pi_\xi)\eta=0$ for
$\eta\in(\mathbb{R}^4)^\perp$. On the other hand, $\mathbb{R}^4$ is invariant under $\mathcal{J}(\pi_\xi)$ and we have already
determined the eigenvalue structure to be $(\mu,0)$ of multiplicities $(2,2)$, respectively, in Lemma \ref{lem-2.3}. Thus
$\mathfrak{C}_\mu^n$ is a complex Osserman model with eigenvalues $(\mu,0)$ of multiplicities $(2,n-2)$, respectively.

Assertions (1) and (2) are invariant under the action of $O(m)$ and by replacing $\xi$ by $J\xi$. Thus, as above, we may assume
$\xi=a_1x_1+b_1Jx_1+b_2Jx_2$ for
$a_1\ne0$. We use the analysis used to prove Lemma \ref{lem-2.3}. We have
$$
J\xi=-b_1x_1-b_2x_2+a_1Jx_1\quad\text{and}\quad L\xi=a_1x_2-b_2Jx_1+b_1Jx_2\,.
$$
We have the following two chains of equivalences:
$$\begin{array}{rrrrrrr}
&\mathcal{J}(\pi_\xi)\xi=\mu\xi&\Leftrightarrow&\xi\perp L\pi_\xi&\Leftrightarrow&J\xi\perp L\xi\phantom{;}\\
\Leftrightarrow&a_1b_2=0&\Leftrightarrow&b_2=0&\Leftrightarrow&\xi\in\pi_{x_1};\\
&\mathcal{J}(\pi_\xi)\xi=0&\Leftrightarrow& \xi\in L\pi_\xi&\Leftrightarrow& J\xi=\pm L\xi\phantom{.}\vphantom{\vrule height 14pt}\\
\Leftrightarrow&b_1=0&\text{and}& a_1=\pm b_2&\Leftrightarrow& \xi=\textstyle\frac1{\sqrt2}(x_1\pm Jx_2).
\end{array}$$
Assertions (1) and (2) now follow.
\end{proof}

\section{Complex Osserman K\"ahler models}\label{sect-3}
In this section, we present some general results we shall need subsequently.
\subsection{Basic results} We refer to \cite{bv-gr-g08} for the proof of the following Lemma:
\begin{lemma}\label{lem-3.1}
Let $\mathfrak{H}$ be a K\"ahler model. Then:
\begin{enumerate}
\item $J^*A=A$.
\item $\mathcal{J}(\pi)$ is complex for all $\pi\in\mathbb{CP}(\mathfrak{H})$.
\item $\mathcal{J}(\pi)y=A(x,Jx)Jy$.
\item If $\mathcal{J}(\pi)=0$ for all $\pi\in\mathbb{CP}(\mathfrak{H})$, then $A=0$.
\end{enumerate}
\end{lemma}

If $\mathfrak{H}$ is complex Osserman K\"ahler and if $\lambda$ is an eigenvalue of $\mathfrak{H}$,
let  $E_\lambda(\pi)$ be the corresponding eigenspace of $\mathcal{J}(\pi)$ for $\pi\in\mathbb{CP}(\mathfrak{H})$.

\begin{lemma}\label{lem-3.2}
Let $\mathfrak{H}$ be a complex Osserman K\"ahler model which has 2 eigenvalues $\lambda_1<\lambda_2$. Then:
\begin{enumerate}
\item $\mathfrak{H}$ is Einstein.
\item If $x,y\in
S(\mathfrak{H})$, then
$y\in E_{\lambda_i}(\pi_x)$ implies $x\in E_{\lambda_i}(\pi_y)$.
\end{enumerate}\end{lemma}

\begin{proof} Let $x\in S(\mathfrak{H})$. By Lemma \ref{lem-3.1}, $J^*A=A$. Thus
$$A(y,Jx,Jx,z)=A(Jy,x,x,Jz)\quad\text{so}\quad\mathcal{J}(Jx)=-J\mathcal{J}(x)J\,.$$
Consequently the Ricci tensor satisfies
$$\rho(x,x)=\operatorname{Tr}\{\mathcal{J}(x)\}=\operatorname{Tr}\{\mathcal{J}(Jx)\}
=\textstyle\frac12\operatorname{Tr}\{\mathcal{J}(\pi_x)\}=\textstyle\frac12
\sum_i\lambda_i\dim\{E_{\lambda_i}(\pi_x)\}\,.$$
As $\rho(x,x)$ is constant on
$S(\mathfrak{H})$,
$\rho(\cdot,\cdot)=c\,\g\cdot\cdot$ so $\mathfrak{H}$ is
Einstein. We suppose $i=2$ so $\lambda_2$ is the maximal eigenvalue; the case $i=1$ is similar. We have
\[
\lambda_2=\max_{z\in S(\mathfrak{H})} \g{\jxa z}{z},
\]
and if $z$ is a unit vector which realizes the maximum, then $z$ is an eigenvector. Hence
Assertion (2) follows from the following sequence of equalities:
\medbreak\qquad\qquad\qquad
$\lambda_2=\g{\jxa y}{y}=A(y,x,x,y)+A(y,Jx,Jx,y)$
\smallbreak\qquad\qquad\qquad
$\phantom{\lambda}=A(x,y,y,x)+A(x,Jy,Jy,x)=\g{\jya x}{x}$.
\end{proof}

\subsection{Eigenvalue multiplicities for complex Osserman models} Methods of algebraic topology can be used to control the
eigenvalue structure of a complex Osserman model. In particular, no
more than $3$ distinct eigenvalues may occur. We refer to \cite{bv-g08} for the proof of the following result:

\begin{theorem}\label{thm-3.3}
If $\mathfrak{H}$ is complex
Osserman, then one of the following holds:
\begin{enumerate}
\item There is just 1 eigenvalue.
\item There are $2$ eigenvalues with multiplicities $(n-2,2)$ with $n\ge4$.
\item There are $2$ eigenvalues with multiplicities $(n-4,4)$ with $n\ge8$.
\item There are $3$ eigenvalues with multiplicities $(n-4,2,2)$ with $n\ge8$.
\end{enumerate}
\end{theorem}

We now impose the K\"ahler identity and apply the relations of
Lemma \ref{lem-3.1}. We begin with a simple observation.
\begin{lemma}\label{lem-3.4}
If $\mathfrak{H}$ is a complex Osserman K{\"a}hler  model with only one eigenvalue, then $\mathfrak{H}$ is flat.
\end{lemma}

\begin{proof}
Suppose that the complex Jacobi operator has only $1$ eigenvalue $\mu$. Then
$$\langle\mathcal{J}(\pi_x) x,x\rangle=A(x,Jx,Jx,x)=s(x)$$ is constant on $S(\mathfrak{H})$ so $\mathfrak{H}$ has  constant
holomorphic sectional curvature. Thus by Lemma~\ref{lem-1.7}, $\mathfrak{H}$ is isomorphic to $\mathfrak{B}_\mu^n$ for some
$\mu$. By Lemma
\ref{lem-2.2}, the eigenvalues of $\mathfrak{B}_\mu^n$ are $\{\frac12\mu,\mu\}$. Thus $\frac12\mu=\mu$ so $\mu=0$ and
$\mathfrak{H}$ is flat.\end{proof}

\begin{remark}\rm
It follows from Equation (\ref{eqn-3}) that $\frac\mu2A_0+\frac\mu6A_J$ is a complex Osserman tensor with constant sectional curvature which
has only one eigenvalue; this does not contradict Lemma \ref{lem-3.4} since this tensor is not K\"ahler.
\end{remark}

\subsection{Critical points of the sectional curvature} We assume that $\mathfrak{H}$ has 2 eigenvalues
henceforth; this is, of course, the case if $n=4$ and $\mathfrak{H}$ is not flat.

\begin{lemma}\label{lem-3.6}
Let $\mathfrak{H}$ be a complex Osserman K{\"a}hler  model which has 2 eigenvalues $\lambda_1<\lambda_2$. Let
$x\in S(\mathfrak{H})$. Then $x$ is a critical point of the holomorphic sectional curvature function
$s$ if and only if $x\in E_{\lambda_i}(x)$ for some $i$.
\end{lemma}

\begin{proof} Suppose $y\in S(E_{\lambda_i}(x))$. Let $\alpha:=\langle y,x\rangle$. Then
$$
R(x,Jx,Jy,x)=R(x,Jx,Jx,y)=\langle\mathcal{J}(\pi_x)y,x\rangle=\lambda_i\alpha\,.
$$
Consider the variation
$v_y(\varepsilon):=(1+2\varepsilon\alpha+\varepsilon^2)^{-1/2}(x+\varepsilon y)\in
S(\mathfrak{H})$. Expand:
\begin{eqnarray*}
&&s(v_y(\varepsilon))=(1+2\varepsilon\alpha)^{-2}\big\{s(x)+
\varepsilon R(x,Jx,Jx,y)+\varepsilon R(x,Jx,Jy,x)\\
&&\qquad\qquad+\varepsilon R(x,Jy,Jx,x)+\varepsilon R(y,Jx,Jx,x)\big\}+O(\varepsilon^2)\\
&&\qquad=(1-4\varepsilon\alpha)\{s(x)+4\varepsilon R(x,Jx,Jx,y)\}+O(\varepsilon^2)\\
&&\qquad=s(x)+\varepsilon\{-4\alpha
s(x)+4\alpha\lambda_i\}+O(\varepsilon^2),\\
&&\partial_{\varepsilon}s(v_y(\varepsilon))|_{\varepsilon=0}=4\alpha(\lambda_i-s(x))\,.
\end{eqnarray*}

Suppose that $x$ is a critical point of the sectional curvature function. Expand
$x=\alpha_1y_1+\alpha_2y_2$ for $y_i\in S(E_{\lambda_i}(x))$. Then
$\partial_{\varepsilon}s(v_{y_i}(\varepsilon))|_{\varepsilon=0}=0$. Since $x$ is
non-zero, at least one of the $\alpha_i$ must be non-zero. Thus
$s(x)=\lambda_i$. Furthermore, if $\lambda_j$ is the other eigenvalue, then $\lambda_i\ne\lambda_j$ so $\lambda_j-s(x)\ne0$ so
$\alpha_j=0$ and thus
$x\in E_{\lambda_i}(x)$.

Conversely, suppose that $x\in E_{\lambda_i}(x)$. If $y\in E_{\lambda_i}(x)$, then $\lambda_i-s(x)=0$
and $\partial_{\varepsilon}(s(v_y(\varepsilon))|_{\varepsilon=0}=0$.
If $y\in E_j(x)$ for $i\ne j$, then $\alpha_j=0$ and so
$\partial_{\varepsilon}(s(v_y(\varepsilon))|_{\varepsilon=0}=0$ as well. Thus we conclude
$\partial_{\varepsilon}(s(v_y(\varepsilon))|_{\varepsilon=0}=0$ for all such variations. Since the
derivatives of all such variations form a spanning set for the tangent space
$T_xS(\mathfrak{H})$, we conclude $x$ is a critical point of the sectional curvature
function.\end{proof}

\begin{lemma}\label{lem-3.7}
Let $\mathfrak{H}$ be a complex Osserman K\"ahler model which has 2
eigenvalues $\lambda_1<\lambda_2$ and which does not have constant holomorphic sectional curvature.
\begin{enumerate}
\item For $i=1,2$, there exist $x_i\in S(\mathfrak{H})$ so that $x_i\in E_{\lambda_i}(x_i)$ and
$s(x_i)=\lambda_i$.
\item If $\dim E_{\lambda_i}\ge4$, then $\lambda_i=0$.
\item We do not have both $\dim E_{\lambda_1}\ge4$ and $\dim E_{\lambda_2}\ge4$.
\end{enumerate}\end{lemma}

\begin{proof} Let $s$ attain its minimum at $x_1\in S(\mathfrak{H})$ and its maximum at $x_2\in S(\mathfrak{H})$.
Since
$\mathfrak{H}$ does not have constant holomorphic sectional curvature, $s(x_1)<s(x_2)$. As $x_1$ and
$x_2$ are critical points of $s$, Assertion (1) follows from Lemma \ref{lem-3.6}; note that $x\in
E_{\lambda_i}(x)$ implies $s(x)=\lambda_i$.

Suppose that $\dim E_{\lambda_i}\ge4$. Choose $x_i$ so $x_i\in E_{\lambda_i}(x_i)$. Since $\dim
E_{\lambda_i}(x_i)\ge4$, there is
$z_i\in E_{\lambda_i}(x_i)$ with
$z_i\perp\pi_{x_i}$. Let $x_i(\theta):=\cos\theta x_i+\sin\theta z_i\in E_{\lambda_i}(x_i)$. By
Lemma \ref{lem-3.2} (2), $x_i\in E_{\lambda}(x_i(\theta))$. We use Lemma \ref{lem-3.1} (3) to see:
\[
\begin{array}{rcl}
\lambda_ix_i&=&A(\cos\theta x_i+\sin \theta z_i,\cos\theta Jx_i+\sin \theta Jz_i)Jx_i\\
\noalign{\smallskip}&=& \cos ^2 \theta A(x_i,Jx_i)Jx_i+2
\cos\theta\sin\theta A(x_i,Jz_i)Jx_i+\sin ^2\theta A(z_i,Jz_i)Jx_i\\
\noalign{\smallskip}&=& \lambda_ix_i+2 \cos\theta\sin\theta A(x_i,Jz_i)Jx_i\,.
\end{array}
\]
Thus $A(x_i,Jz_i)Jx_i=A(z_i,Jx_i)Jx_i=0$ and, similarly, $A(z_i,x_i)x_i=0$. Hence
$$\lambda_iz_i=A(x_i,Jx_i)z_i=A(z_i,x_i)x_i+A(z_i,Jx_i)Jx_i=0\,.$$
This shows $\lambda_i=0$ and establishes Assertion (2); Assertion (3) follows from Assertion (2)
since $0=\lambda_1<\lambda_2=0$ is not possible.
\end{proof}

\subsection{The proof of Theorem \ref{thm-1.3}}
Let $\mathfrak{H}$ be a complex Osserman K\"ahler model which is not
flat. By Lemma
\ref{lem-3.7}, the eigenvalue multiplicity
$(n-4,4)$ with $n\ge8$ is not possible. Lemma \ref{lem-3.4} shows that the multiplicity $(n)$ is
not possible. Theorem~\ref{thm-3.3} then shows the multiplicities to be $(n-2,2)$ or $(n-4,2,2)$
with $n=4k\ge8$.\hfill\qed
\section{The proof of Theorem \ref{thm-1.8}}\label{sect-4}
Through out this section, let  $\mathfrak{H}$ be a complex Osserman K\"ahler model with 2 eigenvalues $(\mu,\lambda)$ of
multiplicities
$(2,n-2)$, respectively, which does not have constant holomorphic sectional curvature. If $n\ge6$, then $\lambda=0$ by
Lemma \ref{lem-3.7} and hence $\mu\ne0$. On the other hand, if $n=4$, we may assume without loss of generality that the notation is
chosen so that
$\mu\ne0$ since both $\mu$ and $\lambda$ can not vanish simultaneously. Thus we shall always assume that $\mu\ne0$ henceforth.

\begin{de}
\rm Complex lines $\{\pi_1,...,\pi_k\}$ in
$\mathbb{CP}(\mathfrak{H})$ will be said to {\it form a $\mu$-configuration} if $\pi_i\subset E_\mu(\pi_i)$ and if
$\pi_i\perp\pi_j$ for $i\ne j$; this then implies $\pi_j\subset E_\lambda(\pi_i)$ for $i\ne j$.\footnote{If
$n=4$, it is a priori possible to have both a $\mu$-configuration and  a $\lambda$-configuration.}
\end{de}

\begin{lemma}\label{lem-4.2} Given $\mathfrak{H}$ as above, there exist complex lines $\{\pi_1,...,\pi_\frac n2\}$ which form
a $\mu$-configuration.
\end{lemma}

\begin{proof} Suppose first $n=4$. Since $\mathfrak{H}$ does not have constant holomorphic sectional curvature, we may apply Lemma \ref{lem-3.7} to
choose
$x\in S(\mathfrak{H})$ so that $x\in E_\mu(\pi_x)$. Let $\pi_1:=\pi_x$ and let $\pi_2:=\pi_1^\perp$. We have
$\pi_2=E_\lambda(\pi_1)$ and hence, dually, $\pi_1=E_\lambda(\pi_2)$ by Lemma \ref{lem-3.2}. Thus $\pi_2=E_\mu(\pi_2)$ and
$\{\pi_1,\pi_2\}$ form a $\mu$-configuration.

Suppose next that $n>4$. We proceed by induction on $n$. By Lemma \ref{lem-3.7}, $\lambda=0$. Use Lemma
\ref{lem-3.7} to choose
$\pi_1$ so
$\pi_1=E_\mu(\pi_1)$. Let
$$\mathfrak{H}_1:=(\pi_1^\perp,\langle\cdot,\cdot\rangle|_{\pi_1^\perp},J|_{\pi_1^\perp},A_{\pi_1^\perp})$$
be the restriction of the model $\mathfrak{H}$ to the subspace $\pi_1^\perp$. As the restriction of a K\"ahler model to a complex subspace
is K\"ahler, we have $\mathfrak{H}_1$ is a K\"ahler model.
If $y\in
S(\pi_1^\perp)$, then $y\in E_0(\pi_1)$ and hence dually $\pi_1\subset E_0(\pi_y)$. This shows $\mathcal{J}(\pi_y)$ preserves $\pi_1$
and hence as $\mathcal{J}(\pi_y)$ is self-adjoint, $\mathcal{J}(\pi_y)$ preserves $\pi_1^\perp$. Furthermore,
$\mathcal{J}(\pi_y)$ has eigenvalues
$(\mu,0)$ of multiplicities $(2,n-4)$ on $\pi_1^\perp$.   Thus  $\mathfrak{H}_1$ is a complex Osserman K\"ahler model of dimension $n-2$ with eigenvalues $(\mu,0)$ of multiplicities $(2,n-4)$. Since
$\lambda=0$ and
$\mu\ne0$, $\mathfrak{H}_1$ does not have constant holomorphic sectional curvature. Consequently we may proceed inductively to
construct a $\mu$-configuration $\{\pi_2,...,\pi_{\frac12n}\}$ for $\mathfrak{H}_1$;
$\{\pi_1,...,\pi_{\frac12n}\}$ is then a $\mu$-configuration for $\mathfrak{H}$.\end{proof}

If $n=4$, then Theorem \ref{thm-1.8} follows from the following result:

\begin{lemma}\label{lem-4.3}
Let $n=4$ and let $\{\pi_1,\pi_2\}$ be a $\mu$-configuration for $\mathfrak{H}$. Fix $x_1\in S(\pi_1)$. There exists $x_2\in
S(\pi_2)$ and constants
$\alpha,\beta$ with
$\alpha+\beta=\lambda$ so that:
\begin{enumerate}
\item $\mathcal{J}(x_1)x_2=\alpha x_2$,\ $\mathcal{J}(Jx_1)x_2=\beta x_2$,\
$\mathcal{J}(x_1)Jx_2=\beta x_2$,\ $\mathcal{J}(Jx_1)Jx_2=\alpha x_2$,\newline
$\mathcal{J}(x_2)x_1=\alpha x_1$,\ $\mathcal{J}(Jx_2)x_1=\beta x_1$,\
$\mathcal{J}(x_2)Jx_1=\beta x_1$,\ $\mathcal{J}(Jx_2)Jx_1=\alpha x_1$.
\item The non-zero curvatures of $A$ are given up to the usual $\mathbb{Z}_2$ symmetries by:
\newline
$A(x_1,Jx_1,Jx_1,x_1)=A(x_2,Jx_2,Jx_2,x_2)=\mu$,\quad$A(x_1,Jx_1,Jx_2,x_2)=\lambda$,
\newline
$A(x_1,x_2,x_2,x_1)=A(Jx_1,Jx_2,Jx_2,Jx_1)=A(Jx_1,Jx_2,x_2,x_1)=\alpha$,\newline
$A(x_1,Jx_2,Jx_2,x_1)=A(Jx_1,x_2,x_2,Jx_1)=A(x_2,Jx_1,Jx_2,x_1)=\beta$.
\item We can choose $x_2\in S(\pi_2)$ so that $\beta=\frac12\mu$ and so that $\alpha=-\frac12\mu$; the model $\mathfrak{H}$
is then equal to $\mathfrak{C}_\mu^4$.
\end{enumerate}
\end{lemma}
\begin{proof} Since $\mathcal{J}(x_1)x_1=0$ and since $\mathcal{J}(x_1)Jx_1=\mathcal{J}(\pi_1)Jx_1=\mu Jx_1$,
$\mathcal{J}(x_1)$ preserves $\pi_1$ and hence as $\mathcal{J}(x_1)$ is self-adjoint,
$\mathcal{J}(x_1)$ preserves $\pi_2$ as well. Thus we may choose an orthonormal basis $\{x_2,Jx_2\}$ for $\pi_2$ so that
$$\mathcal{J}(x_1)x_2=\alpha x_2\quad\text{and}\quad\mathcal{J}(x_1)Jx_2=\beta Jx_2\,.$$
With this choice of basis, we have that
\medbreak\qquad
$A(x_2,x_1,x_1,x_2)=\langle\mathcal{J}(x_1)x_2,x_2\rangle=\alpha$,
\smallbreak\qquad
$A(Jx_2,x_1,x_1,Jx_2)=\langle\mathcal{J}(x_1)Jx_2,Jx_2\rangle=\beta$,
\smallbreak\qquad
$A(Jx_2,x_1,x_1,x_2)=\langle\mathcal{J}(x_1)Jx_2,x_2\rangle=0$.
\medbreak\noindent
Similarly $\mathcal{J}(Jx_1)$ preserves $\pi_2$, $\mathcal{J}(x_2)$ preserves $\pi_1$, and $\mathcal{J}(Jx_2)$ preserves
$\pi_1$. Since
$\{x_i,Jx_i\}$ is an orthonormal basis for $\pi_i$, the remaining identities of Assertion (1) follow from the calculations:
\medbreak\qquad
$\langle\mathcal{J}(Jx_1)x_2,x_2\rangle=A(x_2,Jx_1,Jx_1,x_2)=A(Jx_2,x_1,x_1,Jx_2)=\beta$,
\smallbreak\qquad
$\langle\mathcal{J}(Jx_1)x_2,Jx_2\rangle=A(x_2,Jx_1,Jx_1,Jx_2)=-A(Jx_2,x_1,x_1,x_2)=0$,
\smallbreak\qquad
$\langle\mathcal{J}(Jx_1)Jx_2,Jx_2\rangle=A(Jx_2,Jx_1,Jx_1,Jx_2)=A(x_2,x_1,x_1,x_2)=\alpha$,
\smallbreak\qquad
$\langle\mathcal{J}(x_2)x_1,x_1\rangle=A(x_1,x_2,x_2,x_1)=\alpha$,
\smallbreak\qquad
$\langle\mathcal{J}(x_2)x_1,Jx_1\rangle=A(x_1,x_2,x_2,Jx_1)=-A(x_2,x_1,x_1,Jx_2)=0$,
\smallbreak\qquad
$\langle\mathcal{J}(x_2)Jx_1,Jx_1\rangle=A(x_2,Jx_1,Jx_1,x_2)=A(Jx_2,x_1,x_2,Jx_2)=\beta$,
\smallbreak\qquad
$\langle\mathcal{J}(Jx_2)x_1,x_1\rangle=A(Jx_2,x_1,x_1,Jx_2)=\beta$,
\smallbreak\qquad
$\langle\mathcal{J}(Jx_2)x_1,Jx_1\rangle=A(x_1,Jx_2,Jx_2,Jx_1)= A(x_2,x_1,x_1,Jx_2)=0$,
\smallbreak\qquad
$\langle\mathcal{J}(Jx_2)Jx_1,Jx_1\rangle=A(Jx_2,Jx_1,Jx_1,Jx_2)=A(x_2,x_1,x_2,x_2)=\alpha$.

\medbreak The curvatures listed in Assertion (2) follow from  the facts that $\mathfrak{H}$ is K\"ahler, that $\pi_i\in
E_\mu(\pi_i)$, that $\pi_i\in E_\lambda(\pi_j)$ for $i\ne j$, and from Assertion (1). We consider
possible missing terms. If there is only  one `$x_i$' term, the curvature must have a form like
$R(x_1,Jx_1,x_1,x_2)$. Such a term vanishes since
$\mathcal{J}(x_1)x_2\perp Jx_1$. Thus we may assume there are two $x_1$ terms and two $x_2$ terms. If two $J$ terms touch  in either the first or the last two arguments, we can use the
K\"ahler identity to remove a
$J$. If there is  one $J$ term it must look like $A(x_1,x_2,x_1,Jx_2)$ modulo the K\"ahler identity; this vanishes by Assertion
(1). The terms with no $J$ in them are $A(x_1,x_2,x_1,x_2)$ which is known. Thus there are two $J$ terms $A(J*,*,J*,*)$
which, modulo the K\"ahler identity, are of the form
$A(Jx_1,x_2,Jx_1,x_2)$ or $A(Jx_1,x_1,Jx_2,x_2)$ which have already been discussed. This proves Assertion (2).

Consider $z=(x_1+x_2)/\sqrt2$. We then have
\begin{eqnarray*}
\mathcal{J}(\pi_z)&=&A(z,Jz)J=\textstyle\frac12A(x_1+x_2,Jx_1+Jx_2)J\\
&=&\textstyle\frac12A(x_1,Jx_1)J+\frac12A(x_2,Jx_2)J+A(x_1,Jx_2)J\\
&=&\textstyle\frac12\{\mathcal{J}(\pi_1)+\mathcal{J}(\pi_2)\}+A(x_1,Jx_2)J\,.
\end{eqnarray*}
It is clear that $\mathcal{J}(\pi_1)+\mathcal{J}(\pi_2)=(\lambda+\mu)\operatorname{Id}$. We compute:
\begin{eqnarray*}
&&A(x_1,Jx_2)Jx_1=A(x_2,Jx_1)Jx_1=\mathcal{J}(Jx_1)x_2=\beta x_2,\\
&&A(x_1,Jx_2)Jx_2=\mathcal{J}(Jx_2)x_1=\beta x_1\,.
\end{eqnarray*}
It now follows that
$$
\mathcal{J}(\pi_z)x_1=\textstyle\frac12(\lambda+\mu)x_1+\beta x_2,\quad
\mathcal{J}(\pi_1)x_2=\textstyle\frac12(\lambda+\mu)x_2+\beta x_1\,.
$$
Since $\mathcal{J}(\pi_z)J=J\mathcal{J}(\pi_z)$, we obtain as well
$$
\mathcal{J}(\pi_z)Jx_1=\textstyle\frac12(\lambda+\mu)Jx_1+\beta Jx_2,\quad
\mathcal{J}(\pi_1)Jx_2=\textstyle\frac12(\lambda+\mu)Jx_2+\beta Jx_1\,.
$$
Thus the eigenvalues of $\mathcal{J}(\pi_z)$ are $\textstyle\frac12(\lambda+\mu)\pm\beta$.
We interchange the roles of $\alpha$ and $\beta$ by considering instead $\tilde z:=(x_1+Jx_2)/\sqrt2$. A similar argument shows the
eigenvalues of
$\mathcal{J}(\pi_{\tilde z})$ are
$\textstyle\frac12(\lambda+\mu)\pm\alpha$. It now follows that
$$\alpha=\pm\beta\,.$$

Suppose first that $\alpha=\beta$. Then $\lambda=\alpha+\beta=2\alpha$ so $\alpha=\beta=\frac12\lambda$. Thus the eigenvalues of
$\mathcal{J}(\pi_z)$ are $\frac12(\lambda+\mu)-\frac12\lambda=\frac12\mu$ and $\frac12(\lambda+\mu)+\frac12\lambda=\lambda+\frac12\mu$. Since
$\mu\ne0$, $\lambda+\frac12\mu\ne\lambda$ and thus we must have $\lambda+\frac12\mu=\mu$ so $\lambda=\frac12\mu\ne0$. We substitute the
values
$\alpha=\beta=\frac\lambda2=\frac\mu4$ into Assertion (2) to determine the curvature  tensor and see it agrees with the expression
given in Example \ref{exm-1.5}. This is not possible as we assumed $\mathfrak{H}$ does not have constant holomorphic
sectional curvature.

Next suppose that $\alpha=-\beta$. We then have $\lambda=0$. By interchanging the roles of $x_2$ and $Jx_2$ if need be, we may
assume that
$\frac12\mu+\beta=\mu$ and $\frac12\mu-\beta=0$. We substitute the values $-\alpha=\beta=\frac12\mu$ into Assertion (2) to determine
the curvature tensor and see it agrees with the expression given in Example
\ref{exm-1.6}.
\end{proof}

We suppose $n\ge6$  henceforth; it then follows from Lemma \ref{lem-3.7} that $\lambda=0$.

\begin{lemma}\label{lem-4.4}
Let $4\le 2k <n$. Let $\{\pi_1,...,\pi_k\}$ be a $\mu$-configuration for $\mathfrak{H}$. Let
$D:=\pi_1+...+\pi_k$ and let
$\mathfrak{D}=(D,\langle\cdot,\cdot\rangle|_D,J|_D,A_D)$.
\begin{enumerate}
\item $\mathfrak{D}$ is a complex Osserman K\"ahler model with 2 eigenvalues $(\mu,0)$ of
multiplicities $(2,2k-2)$, respectively, which does not have constant holomorphic sectional curvature.
\item If $u,v,w\in D$, then $A(u,v)w\in D$ and $\mathcal{J}(u)v\in D$.
\end{enumerate}
\end{lemma}

\begin{proof} We apply the argument used in the proof of Lemma \ref{lem-4.2} to establish Assertion (1). Let
$y\in D^\perp$. Then
$y\in E_0(\pi_i)$ so dually
$\pi_i\subset E_0(\pi_y)$. Thus $D\subset E_0(\pi_y)$. Let $z\in S(D)$. Then $z\in E_0(\pi_y)$ and hence $y\in E_0(\pi_z)$. Thus
$D^\perp\subset E_0(\pi_z)$. Consequently $\mathcal{J}(\pi_z)D^\perp\subset D^\perp$ and as $\mathcal{J}(\pi_z)$ is self
adjoint,
$\mathcal{J}(\pi_z)D\subset D$. It now follows that $\mathfrak{D}$ is a complex Osserman K\"ahler model with eigenvalues
$(\mu,0)$ of multiplicities
$(2,2k-2)$, respectively;
$\mathfrak{D}$ does not have constant holomorphic sectional curvature since $\mu\ne0$ and $\lambda=0$. This proves Assertion (1).

Let $\{x_i,Jx_i\}$ be an orthonormal basis for $\pi_i$; $\{x_1,Jx_1,...,x_k,Jx_k\}$ is an orthonormal basis for $D$. Let $i\ne j$. Let
$z:=\cos(\theta)x_i+\sin(\theta)x_j\in S(D)$. We noted above that $\mathcal{J}(\pi_z)D\subset D$. We expand:
\begin{eqnarray*}
\mathcal{J}(\pi_z)&=&A(\cos(\theta)x_i+\sin(\theta)x_j,\cos(\theta)Jx_i+\sin(\theta)Jx_j)J\\
&=&\cos^2(\theta)\mathcal{J}(\pi_i)+\sin^2\theta\mathcal{J}(\pi_j)+2\cos\theta\sin\theta A(x_i,Jx_j)J\,.
\end{eqnarray*}
This shows that $A(x_i,Jx_j)$ preserves $D$; similarly $A(x_i,x_j)$ and $A(Jx_i,Jx_j)$ preserve $D$. Since $A(x_i,Jx_i)$ also preserves
$D$, we conclude $A(u,v)w\in D$ for any $u,v,w\in D$. Taking $v=w$ shows $\mathcal{J}(u)$ preserves
$D$.
\end{proof}

\medbreak\noindent{\bf The proof of Theorem \ref{thm-1.8}}
Let $n\ge6$. Let $\{\pi_1,...,\pi_{\frac n2}\}$ be a
$\mu$-configuration. Choose $x_1\in S(\pi_1)$. We apply Lemma \ref{lem-4.3} to choose $x_i\in S(\pi_i)$ for $i\ge2$ so that the
non-zero curvatures involving the indices $x_1$ and $x_j$ for $2\le j\le \frac n2$ are given by Example \ref{exm-1.6}.
Lemma \ref{lem-4.4} shows that $A(\cdot,\cdot,\cdot,\cdot)$ vanishes if there are 3 impacted indices. Thus to show
$\mathfrak{H}$ is isomorphic to
$\mathfrak{C}_\mu^n$, we need only examine the curvature tensor on $\pi_j+\pi_k$ for $2\le j<k\le\frac n2$. To simplify the notation,
we set $j=2$ and $k=3$. Let $\mathfrak{D}$ be the K\"ahler model
determined by the
$\mu$-configuration
$\{\pi_1,\pi_2,\pi_3\}$; let $\mathfrak{D}_{ij}$ be the K\"ahler model determined by the $\mu$-configuration
$\{\pi_i,\pi_j\}$ where $i<j$;
$\mathfrak{D}$ and $\mathfrak{D}_{ij}$ are complex Osserman K\"ahler models with eigenvalues
$\{\mu,0\}$ with multiplicities
$(2,4)$ and $(2,2)$, respectively which do not have constant holomorphic sectional curvature. Let
$$
\pi:=\pi_{\frac1{\sqrt3}(x_1+x_2+x_3)},\quad \pi_{ij}:=\pi_{\frac1{\sqrt2}(x_i+x_j)}\quad\text{for}\quad i<j\,.
$$
We may expand
\begin{eqnarray}
&&\mathcal{J}(\pi_i)=A(x_i,Jx_i)J,\nonumber\\
&&2\mathcal{J}(\pi_{ij})=A(x_i,Jx_i)J+A(x_j,Jx_j)J+2A(x_i,Jx_j)J,\label{eqn-5}\\
&&3\mathcal{J}(\pi)=\textstyle\sum_{1\le i\le 3}A(x_i,Jx_i)J+2\sum_{1\le i<j\le 3}A(x_i,Jx_j)J\,.\nonumber
\end{eqnarray}
Since $\mathcal{J}(\pi_1)+\mathcal{J}(\pi_2)+\mathcal{J}(\pi_3)=\mu\operatorname{Id}$ on $D$, we may express
\begin{equation}\label{eqn-6}
3\mathcal{J}(\pi)=2\mathcal{J}(\pi_{12})+2\mathcal{J}(\pi_{13})+2\mathcal{J}(\pi_{23})-\mu\operatorname{Id}\,.
\end{equation}
By replacing $A$ by $-A$ if necessary, we may assume without loss of generality that $\mu>0$. As $\mu$ is the largest
eigenvalue, if $\sigma$ is any complex line, $\langle\mathcal{J}(\sigma)w,w\rangle\le\mu|w|^2$ and equality holds if and only
if
$\mathcal{J}(\sigma)w=\mu w$. Choose $w\in S(D)$ so $\mathcal{J}(\pi)w=\mu w$. By replacing $w$ by $\cos\theta
w+\sin\theta Jw$ if necessary, we may suppose $\langle w,Jx_1\rangle=0$ so
$w=a_1x_1+a_2x_2+a_3x_3+b_2Jx_2+b_3Jx_3$.
Set
\begin{eqnarray*}
&&w_{12}:=a_1x_1+a_2x_2+b_2Jx_2,\quad w_{13}:=a_1x_1+a_3x_3+b_3Jx_3,\\
&&w_{23}:=a_2x_2+a_3x_3+b_2Jx_2+b_3Jx_3\,.
\end{eqnarray*}
We then have $|w_{12}|^2+|w_{13}|^2+|w_{23}|^2=2$. Since the curvature vanishes if there are 3 impacted indices,
$\mathcal{J}(\pi_{ij}) w=\mathcal{J}(\pi_{ij})w_{ij}$. We use Equation
(\ref{eqn-6}) to estimate:
\begin{eqnarray*}
4\mu&=&\langle(3\mathcal{J}(\pi)+\mu\operatorname{Id})w,w\rangle
   =2\textstyle\sum_{i<j}\langle\mathcal{J}(\pi_{ij})w_{ij},w_{ij}\rangle\\
&\le&2\mu\{|w_{12}|^2+|w_{13}|^2+|w_{23}|^2\}=4\mu\,.
\end{eqnarray*}
This shows that all the inequalities must in fact have been equalities and thus
$$w_{ij}\in E_\mu(\pi_{ij})\quad\text{for}\quad 1\le i<j\le 3\,.$$

By Lemma \ref{lem-4.3}, the curvature of the models $\mathfrak{D}_{12}$ and $\mathfrak{D}_{13}$ is given by Example
\ref{exm-1.6}. Thus Lemma
\ref{lem-2.5} applies. Since $\pi_{1i}$ is determined by $(x_1+x_i)/\sqrt{2}$, $\pi_{1i}=E_\mu(\pi_{1i})$. Thus
$w_{1i}\in\pi_{1i}$ for $i=2,3$. Since
$b_1=0$, this implies
$a_1=a_2=a_3$ and $b_2=b_3=0$. Consequently $w=\pm(x_1+x_2+x_3)/\sqrt{3}$ and we conclude
$\mathcal{J}(\pi_{23})(x_2+x_3)=\mu(x_2+x_3)$ or, equivalently by Equation (\ref{eqn-5}), that
$$\mu(x_2+x_3)=2A(x_2,Jx_3)J(x_2+x_3)\,.$$
We take the inner product with $x_2$ to conclude $A(x_2,Jx_3,Jx_3,x_2)=\frac12\mu$, i.e.  $\mathcal{J}(x_2)Jx_3=\frac\mu2 Jx_3$. By
Lemma \ref{lem-4.3}, the eigenvalues of $\mathcal{J}(x_2)$ on $\pi_3$ are $\pm\frac\mu2$. Consequently
$\mathcal{J}(x_2)x_3=-\frac\mu2 x_3$. Lemma \ref{lem-4.3} now shows that the curvatures of $\mathfrak{D}_{23}$ are given by
Example \ref{exm-1.6}. This completes the proof.\hfill\qed

\section{The proof of Theorem \ref{thm-1.9}}\label{sect-5}

\begin{de}\rm We say that $(M,g)$ is a {\it regular curvature homogeneous manifold} if given any point $P\in M$, there
exists an open neighborhood $\mathcal{O}$ of $P$ and an or\-tho\-nor\-mal frame field $\{e_1,..., e_n\}$ so that the curvature
tensor $R(e_i,e_j,e_k,e_l)(x)=C_{ijkl}$ is independent of the point $x\in\mathcal{O}$.
\end{de}

\begin{lemma}\label{lem-5.3}
Let $\mathcal{H}=(M,g,J)$ be a complex Osserman manifold whose curvature is modeled on
$\mathfrak{C}_\mu^4$ at every point of $M$. Then $(M,g)$ is a regular curvature homogeneous Einstein manifold.
\end{lemma}

\begin{proof} Fix $P\in M$ and let $\mathcal{O}$ be a contractible open
neighborhood of $P$. Let $R_L:=3\{\textstyle\frac{\mu}{2}A_g+\frac{\mu}{6}A_J-R\}$ define  a smooth algebraic curvature tensor
on $T(\mathcal{O})$. The analysis of Lemma \ref{lem-2.3} shows that $R_L$ is Osserman with eigenvalues $\{0,\mu\}$ of
multiplicities
 $(3,1)$, respectively. Let $\sigma(s)$ be orthogonal projection on $E_\mu(\mathcal{J}_{R_L}(s))$ for $s\in S(TM)$. The Cauchy
integral formula
$$\sigma(s)=\frac1{2\pi i}\int_{|z-\mu|=\frac\mu2}\{\mathcal{J}_{R_L}(s)-z\operatorname{Id}\}^{-1}dz$$
shows that $\sigma$ varies smoothly. Consequently,
$\Sigma:=\operatorname{Range}(\sigma)$ is a smooth line bundle over $S(T\mathcal{O})=\mathcal{O}\times S^3$.
Consequently $\Sigma$ is trivial and admits a smooth unit section
$\mathcal{L}$. Since $\mathcal{L}x=\pm Lx$, this analysis shows that the auxiliary almost complex structure  $L$ of Lemma
\ref{lem-2.3} can be chosen to vary smoothly with the point of the manifold.

We  have $(JL)^2=1$ and $(JL)^*=LJ=JL$ is self-adjoint. Let $E_{\pm1}(JL)$ be the associated $\pm1$ eigenbundles of
$JL$. These comprise orthogonal $J$-invariant $2$-planes. Let
$s$ be a smooth unit tangent vector such that
$Js\ne\pm Ls$. Form:
$$s_\pm:=|s\pm JLs|^{-1}(s\pm JLs)\in S(E_{\pm1}(JL))\quad\text{and}\quad x_1=(s_++s_-)/\sqrt2\,.$$
We then have
 $$g(Jx_1,Lx_1)=-g(JLx_1,x_1)=-\textstyle\frac12 g(s_+-s_-,s_++s_-)=0\,.$$
Let $\pi_1:=\operatorname{Span}\{x_1,Jx_1\}$. Then $\pi_1=E_\mu(\pi_1)$. We use the Cauchy integral formula to choose a smooth
unit section to the $-\frac\mu2$ eigenspace of $\mathcal{J}(x_1)$ on $\pi_1^\perp$. The proof that the curvature
tensor has the requisite form relative to this local frame now follows using the argument to establish Lemma \ref{lem-4.3}.
\end{proof}

Let $(M,g)$ be an Einstein regular curvature homogeneous manifold of dimension  $n=4$.
Let $R_\Lambda$ be the induced action of the curvature tensor on $\Lambda^2$.
Choose a local orientation for
$M$ and decompose
$\Lambda^2=\Lambda^2_+\oplus\Lambda^2_-$. Let $W$ be the Weyl conformal curvature tensor.  Since $(M,g)$ is Einstein,
$R=W+cA_0$ and hence $R_\Lambda=W_\Lambda+c\operatorname{Id}$ for some suitably chosen constant $c$. Since
$W_\Lambda:\Lambda^2_\pm\rightarrow\Lambda^2_\pm$,  $R=R_+\oplus R_-$ where $R_\pm\in\Lambda^2_\pm\otimes\Lambda^2_\pm$.
Let
$\{\lambda_1^\pm,\lambda_2^\pm,\lambda_3^\pm\}$ be the eigenvalues of $R_{\Lambda_\pm}$; these are constant as $(M,g)$ is
curvature homogeneous. We let $\omega_i^\pm$ be an orthonormal frame for $\{\Lambda^2\}$ so that $|\omega_i^\pm|^2=2$ and
$R_{\Lambda_\pm}\omega_i^\pm=\lambda_i^\pm \omega_i^\pm$ where
$$\omega_i^\pm=\textstyle\sum_{ab} \omega^\pm_{i,ab}e_a\wedge e_b\,.$$
Since $(M,g)$ is curvature homogeneous, the coefficients $\omega_{ij,ab}$ are constant. We
specialize an argument of Derdzi\'nski \cite{d83} to show:
\begin{lemma}\label{lem-5.2}
Let $(M,g)$ be an Einstein regular curvature homogeneous manifold of dimension  $n=4$. If
$\lambda_1^+\ne\lambda_2^+=\lambda_3^+$ and if the matrix $\omega_{1,ab}^+$ is invertible, then $\nabla\omega_1^+=0$.
\end{lemma}

\begin{proof} As $(M,g)$ is Einstein, the Ricci tensor $\rho$ is parallel. The second Bianchi identity yields:
\begin{eqnarray*}
(\delta_1R)(y,z,w)&:=&R(e_i,y,z,w;e_i)=-R(e_i,y,w,e_i;z)-R(e_i,y,e_i,z;w)\\
&=&-\rho(y,w;z)+\rho(y,z;w)=0\,.
\end{eqnarray*}
As the decomposition $\Lambda^2=\Lambda_+^2\oplus\Lambda_-^2$ is parallel, we have $\delta_1R_\pm=0$ individually.
Furthermore, there exist smooth $1$-forms $\phi_{ij}^+=-\phi_{ji}^+$ so that
$$
\nabla\omega_i^+=\textstyle\sum_j\phi_{ij}^+\otimes\omega_j^+\quad\text{where}\quad\phi_{ij}^+=\sum_a\phi_{ij,a}e_a\,.
$$
We use the eigenvalue decomposition to express:
$$
R_+=\textstyle\frac12\{\lambda_1^+\omega_1^+\otimes\omega_1^+
+\lambda_2^+\omega_2^+\otimes\omega_2^++\lambda_3^+\omega_3\otimes\omega_3^+\}\,.
$$
Since the eigenvalues $\lambda_i^+$ are constant and since $\phi_{ij}^+=-\phi_{ji}^+$, we have:
\begin{eqnarray*}
\nabla R_+&=&
\textstyle\frac12\sum_{ij}\lambda_i^+\phi_{ij}^+\otimes\{\omega_j^+\otimes\omega_i^++\omega_i^+\otimes\omega_j^+\}\\
&=&\textstyle\frac12\sum_{ij}\{\lambda_i^+-\lambda_j^+\}\phi_{ij}^+\otimes\omega_j^+\otimes\omega_i^+\,,\nonumber\\
\delta_1 R_+&=&\textstyle\frac12\sum_{ijab}\{\lambda_i^+-\lambda_j^+\}
\phi_{ij,a}^+\omega_{j,ab}^+e_b\otimes\omega_i^+\,.\nonumber
\end{eqnarray*}
We set $i=2$ and $i=3$ in this identity. Since $\lambda_2^+=\lambda_3^+$, we set $j=1$ to show:
\begin{eqnarray*}
0&=&\textstyle\frac12\{\lambda_2^+-\lambda_1^+\}\sum_{ab}\phi_{21,a}^+\omega_{1,ab}^+e_b\,,\\
0&=&\textstyle\frac12\{\lambda_3^+-\lambda_1^+\}\sum_{ab}\phi_{31,a}^+\omega_{1,ab}^+e_b\,.
\end{eqnarray*}
As $\omega_{1,ab}^+$ is invertible, $\phi_{21}^+=0$ and $\phi_{31}^+=0$ so
$\nabla\omega_1^+=0$.
\end{proof}

We may combine Lemma \ref{lem-5.2} and Lemma \ref{lem-5.3} to show:

\begin{lemma}\label{lem-5.4}
Let $\mathcal{H}=(M,g,J)$ be a connected $4$-dimensional complex Osserman manifold whose curvature is modeled on
$\mathfrak{C}_\mu^4$ at every point. Then $\nabla R=0$.
\end{lemma}

\begin{proof} Fix $P\in M$. We may use Lemma \ref{lem-5.3} to choose a local orthonormal frame $\{x_1,Jx_1,...\}$ realizing the
curvature relations of Example \ref{exm-1.6}. Let $L$ be defined as above. Define:
\begin{eqnarray*}
&&\phi_1^\pm:=x_1\wedge Jx_1\pm x_2\wedge Jx_2,\qquad
  \phi_2^\pm:=x_1\wedge x_2\mp Jx_1\wedge Jx_2,\\
&&\phi_3^\pm:=x_1\wedge Jx_2\pm Jx_1\wedge x_2\,.
\end{eqnarray*}
This forms an orthogonal frame for $\Lambda_\pm^2$ where $||\phi_i^\pm||^2=2$.
We use the structure equations of Example \ref{exm-1.6} to compute:
\medbreak
$R\phi_1^\pm=\frac12\{R(x_1,Jx_1,x_1,Jx_1)\pm 2 R(x_1,Jx_1,x_2,Jx_2)+R(x_2,Jx_2,x_2,Jx_2)\}\phi_1^\pm$
\smallbreak\quad\hphantom{A.}
$=\frac12\{-\mu\pm0-\mu\}\phi_1^\pm$,
\medbreak
$R\phi_2^\pm=\frac12\{R(x_1,x_2,x_1,x_2)\mp2R(x_1,x_2,Jx_1,Jx_2)+R(Jx_1,Jx_2,Jx_1,Jx_2\}\phi_2^\pm$
\smallbreak\quad\hphantom{A.}
$=\frac12\{\frac\mu2\mp2\frac\mu2+\frac\mu2\}\phi_2^\pm$,
\medbreak
$R\phi_3^\pm=\frac12\{R(x_1,Jx_2,x_1,Jx_2)\pm 2R(x_1,Jx_2,Jx_1,x_2)+R(Jx_1,x_2,Jx_1,x_2)\}\phi_3^\pm$
\smallbreak\quad\hphantom{A.}
$=\frac12\{-\frac\mu2\pm 2\frac\mu2-\frac\mu2\}\phi_3^\pm$.
\medbreak\noindent This yields
$$\begin{array}{rrr}
R\phi_1^+=-\mu\phi_1^+,&R\phi_2^+=0\phi_2^+,&R\phi_3^+=\phantom{-}0\phi_3^+,\\
R\phi_1^-=-\mu\phi_1^-,&R\phi_2^-=\mu\phi_2^-,&R\phi_3^-=-\mu\phi_3^-.
\end{array}
$$
Thus these are eigenvectors of $R$. Since $\phi_1^+$ is  the K\"ahler form of
$J$, $\phi_{1,ab}^+$ is an invertible matrix. As $\lambda_1^+\ne\lambda_2^+=\lambda_3^+$,  Lemma
\ref{lem-5.2} yields $\nabla\phi_1^+=0$ and, equivalently, $\nabla J=0$. Since $\lambda_1^-=\lambda_3^-\ne\lambda_2^-$, and
since the K\"ahler form of
$L$ is
$\phi_2^+$, a similar argument shows $\nabla L=0$. We have $\nabla A_0=0$. Since $\nabla J=0$ and $\nabla L=0$, we have
$\nabla A_J=0$ and $\nabla A_L=0$ as well. It now follows that $\nabla R=0$ at $P$.
Since $P$ was arbitrary, the Lemma follows.
\end{proof}

{\bf The proof of Theorem \ref{thm-1.9}} Suppose $\mathcal{H}=(M,g,J)$ is a connected complex Osserman
manifold which does not have constant holomorphic sectional
curvature at at least one point $Q$ of the manifold; we argue for a contradiction. We
use Theorem \ref{thm-1.8} to see that  at any point $P$ of $M$, the curvature is either modeled on $\mathfrak{B}_\mu^4$ or on
$\mathfrak{C}_\mu^4$. If the curvature at $P$ is modeled on $\mathfrak{B}_\mu^4$, then $\mu=2\lambda\ne0$ while if the
curvature is modeled on $\mathfrak{C}_\mu^4$, then $\lambda=0$. If
$\lambda\ne0$ at some point, then $\lambda\ne0$ on an open set $\mathcal{O}$ and hence $\mu=2\lambda$ on $\mathcal{O}$. Since
$(M,g)$ is Einstein, the scalar curvature is constant. Since the scalar curvature in this setting is $6\mu$, $\mu$ and hence
$\lambda$ are constant on
$\mathcal{O}$ so $\lambda\ne0$ on the closure of $\mathcal{O}$. This implies $\mathcal{O}=M$ and contradicts the assumption
that $\mathcal{H}$ does not have constant holomorphic sectional curvature at some point. Thus at every point, the curvature is
modeled on $C_\mu^4$; since $\tau=4\mu$ we have $\mu$ is constant. We may therefore use Lemma \ref{lem-5.4} to see $\nabla
R=0$. Let $R(u_1,u_2,u_3,u_4;u_5,u_6)$ be the components of $\nabla^2R$. The curvature tensor of a  locally symmetric space is
very restrictive. In particular, we have:
\begin{eqnarray*}
0&=&R(u_1,u_2,u_3,u_4;u_5,u_6)- R(u_1,u_2,u_3,u_4;u_6,u_5)\\
&=&R(R(u_6,u_5)u_1,u_2,u_3,u_4)+R(u_1,R(u_6,u_5)u_2,u_3,u_4)\\
&+&R(u_1,u_2,R(u_6,u_5)u_3,u_4)+R(u_1,u_2,u_3,R(u_6,u_5)u_4)\,.
\end{eqnarray*}
We apply this identity with $u_1=x_1$, $u_2=x_2$, $u_3=x_2$, $u_4=Jx_1$,  $u_5=x_1$, $u_6=Jx_1$ to compute:
\begin{eqnarray*}
0&=&R(R(Jx_1,x_1)x_1,x_2,x_2,Jx_1)+R(x_1,R(Jx_1,x_1)x_2,x_2,Jx_1)\\
&+&R(x_1,x_2,R(Jx_1,x_1)x_2,Jx_1)+R(x_1,x_2,x_2,R(Jx_1,x_1)Jx_1)\,.
\end{eqnarray*}
We substitute the relations of Example \ref{exm-1.6} to conclude:
\begin{eqnarray*}
0&=&\mu R(Jx_1,x_2,x_2,Jx_1)+0+0-\mu R(x_1,x_2,x_2,x_1)\\
&=&\textstyle\frac12\mu^2+\frac12\mu^2\,.
\end{eqnarray*}
This shows that $\mu=0$. Hence $(M,g)$ is flat, contrary to our assumption.\hfill\qed

\section*{Acknowledgements}
The authors would like to thank Prof. E. Garc\'ia-R\'io for many useful conversations.
Research of M. Brozos-V\'azquez supported by projects MTM2009-07756 and INCITE09 207 151 PR (Spain).
Research of P. Gilkey partially supported by project MTM2009-07756 (Spain).

\end{document}